\documentclass{amsproc}

\newtheorem{theorem}{Theorem}[section]
\newtheorem{lemma}[theorem]{Lemma}

\theoremstyle{definition}

\newtheorem{assumption}[theorem]{Assumption}

\theoremstyle{remark}

\numberwithin{equation}{section}

\newcommand{\abs}[1]{\lvert#1\rvert}

\newcommand{\D}{{\rm Dom}}
\newcommand{\R}{{\rm Range}}
\newcommand{\s}{{\rm supp}}

\begin{document}

\title{The unique ergodicity of equicontinuous laminations}

\author{Shigenori Matsumoto}
\address{Department of Mathematics, College of
Science and Technology, Nihon University, 1-8-14 Kanda, Surugadai,
Chiyoda-ku, Tokyo, 101-8308 Japan
}
\curraddr{Department of Mathematics, College of
Science and Technology, Nihon University, 1-8-14 Kanda, Surugadai,
Chiyoda-ku, Tokyo, 101-8308 Japan}
\email{matsumo@math.cst.nihon-u.ac.jp
}
\thanks{The author is partially supported by Grant-in-Aid for
Scientific Research (C) No.\ 20540096.}
\subjclass{Primary 53C12,
secondary 37C85.}

\keywords{lamination, foliation, transversely invariant measure,
unique
ergodicity}

\date{October 7, 2009 }

\begin{abstract}
We prove that a transversely equicontinuous
minimal lamination on a locally compact metric space $Z$
has a transversely invariant Radon measure. Moreover if
the space $Z$ is compact, then the tranversely invariant
Radon measure is shown to be unique up to a scaling.

\end{abstract}

\maketitle

\section{Introduction}

Let $Z$ be a locally compact metric space, $\mathcal L$ a
$p$-dimensional
lamination on $Z$. {\em We assume throughout that $\mathcal L$ is
minimal}. Let $h: \mathbb R^p\times X\to Z$ be a lamination chart,
i.\ e.\ a homeomorphism onto an open subset $h(\mathbb R^p\times X)$
such that the plaque $h(\mathbb R^p\times\{x\})$
lies on a leaf of $\mathcal L$ for any $x\in X$. We identify $X$ with
the image $h(\{0\}\times X)$ and call it a {\em cross section} of 
$\mathcal L$. With the metric induced from $Z$, $X$ is also locally
compact. Notice that any leaf of $\mathcal L$ intersects $X$.

  Given a leafwise curve joining two points $x$ and
$y$
on $X$, a holonomy map along $c$ is defined as usual
to be a local homeomorphism $\gamma$ from an
open neighbourhood $\D(\gamma)$ of $x$ onto an open neighbourhood
$\R(\gamma)$ of $y$. We say that $\mathcal L$ is {\em transversely
equicontinuous} w.\ r.\ t.\ a cross section $X$ if the
family of all the corresponding holonomy maps is equicontinuous.

\begin{theorem} \label{t1}
Let $\mathcal L$ be a minimal lamination
on a locally compact metric space $Z$, transversely
 equicontinuous
w.\ r.\ t.\ a cross section $X$. Then there is a Radon measure
on $X$  which is left invariant by any holonomy map. If further $Z$  is compact, then the invariant measure is
unique up to a scaling.
\end{theorem}

The existence of invariant measure was already shown by R. 
Sackesteder in \cite{S} for a pseudogroup acting on a compact metric
space. But the compactness condition for a cross section
is  too strong 
to obtain a corresponding result for laminations or foliations
(even on compact spaces or manifolds).
In section 2, we include a slightly general theorem applicable to laminations;
the proof closely follows an argument in Lemme 4.4 in \cite{DKN},
which is meant for codimension one 
foliations.

In section 3 we show the uniqueness for a compact lamination. The
argument here which is adapted for pseudo*groups as defined in section
2
 is rather messy, but the
original idea is quite simple, which the reader can find in section 4.
 
In section 4, we deal with an equicontinuous group action on a compact
metric
space, together with a random walk on a group. We show that the
corresponding harmonic probability measure on the space is unique.

The uniqueness of harmonic measures for tangentially sufficiently
smooth foliations and laminations (\cite{C},\cite{G}) remains an open question.

\section{The existence}

Let $Y$ be a Hausdorff space. By a local homeomorphism, we
mean a homeomorphism $\gamma$ from an open subset $\D(\gamma)$
of $Y$ onto an open subset $\R(\gamma)$.
A set $\Gamma$ of local homeomophisms of $Y$ is called a {\em pseudo*group},
if it satisfies the following conditions.
\\
(1) If $\gamma\in\Gamma$ and $U$ is an open subset of $\D(\gamma)$,
then the restriction $\gamma\vert_U$ is in $\Gamma$.
\\
(2) The identity $id_X$ belongs to $\Gamma$.
\\
(3) If $\gamma,\gamma'\in\Gamma$ and $\D(\gamma')=\R(\gamma)$,
then the composite $\gamma'\circ\gamma$ is in $\Gamma$.
\\
(3) If $\gamma\in\Gamma$, then $\gamma^{-1}\in\Gamma$.

This differs from the usual definition of pseudogroups in
that it does not assume the axiom for taking the union.
Thus for example the set of all the holonomy maps w.\ r.\ t.\
a cross section given in section 1 forms a pseudo*group,
while the pseudogroup they generate might be bigger.
There are two reasons for introducing the concepts of pseudo*groups:
one is that in Theorem \ref{t1}, assuming the equicontinuity
for the pseudogroup generated by the holonomy maps may be
stronger than what we have tacitly in mind: the other is that
some part of the argument in section 3 cannot be put into
the framework of the usual pseudogroups.

\bigskip

Let $X$ be a locally compact metric space and $\Gamma$ a pseudo*group
of local homeomorphsims of $X$. We assume that
the action is minimal, i.\ e.\ the $\Gamma$-orbit of any point is dense
in $X$, and that the action is equicontinuous, i.\ e.\ for any
$\epsilon>0$, there is $\delta(\epsilon)>0$ such that if
$\gamma\in\Gamma$, $x,x'\in\D(\gamma)$ and $d(x,x')<
\delta(\epsilon)$,
then we have $d(\gamma x,\gamma x')\leq \epsilon$.

Denote by $C_c(X)$ the space of real valued continuous functions $\zeta$
whose support $\s\zeta$ is compact.
 A Radon measure $\mu$ on $X$
is called $\Gamma$-invariant if whenever $\zeta\in C_c(X)$ and
$\gamma\in\Gamma$ satisfy $\s\zeta\subset\D(\gamma)$, we have
$\mu(\zeta\circ\gamma^{-1})=\mu(\zeta)$.
In fact if $\mu$ is $\Gamma$-invariant, we get a bit more, e.\ g.\
for any bounded continuous function $\zeta: X\to{\mathbb R}$ which
vanishes outside $\D(\gamma)$, we have
$\mu(\zeta\circ\gamma^{-1})=\mu(\zeta)$,
as the dominated convergence theorem shows. In this case the both
hand sides might be $\infty$. This will be used 
in section 3.

 Let $X_0$ be
a relatively compact open subset of $X$, and denote by $\Gamma_0$
the restriction of $\Gamma$ to $X_0$ i.\ e.\
$$
\Gamma_0=\{\gamma\in\Gamma\mid \D(\gamma)\cup \R(\gamma)\subset X_0\}.$$

The purpose of this
section is to show the following theorem.

\begin{theorem} \label{t2}
There exists a finite $\Gamma_0$-invariant Radon measure $\mu$ on $X_0$.
\end{theorem}

The minimality assumption shows then the existence of $\Gamma$-invariant
measure on $X$ and the proof of the existence part of Theorem \ref{t1}
will be complete.

\newcommand{\nn}{C_c(X)_{\geq0}}
\newcommand{\nnn}{C_c(X_0)_{\geq0}}
\newcommand{\p}{C_c(X)_{>0}}
\newcommand{\pp}{C_c(X_0)_{>0}}

\bigskip

Let us define
$$
\nn=\{\zeta\in C_c(X)\mid \zeta\geq 0\}\ \ \ \rm{and}
$$
$$
\p=\{\zeta\in \nn\mid \zeta(x)>0, \ \ \exists x\in X\}.
$$

For any $\psi\in C_c(X)$ and $\gamma\in\Gamma$, extend the function
$\psi\circ\gamma^{-1}$ to the whole $X$ so as to vanish outside
$\R(\gamma)$ and still denote it by $\psi\circ\gamma^{-1}$. It may
no longer be continuous.
For any $\zeta\in\nn$ and $\psi\in\p$, define $(\zeta:\psi)$ by
$$
(\zeta:\psi)=\inf\{\sum_{i=1}^nc_i\mid\zeta\leq\sum_{i=1}^nc_i
\psi\circ\gamma_i^{-1},\ c_i>0,\ \gamma_i\in\Gamma,\ n\in{\mathbb N} \}.
$$
Notice that the minimality of $\Gamma$ implies that
$(\zeta:\psi)<\infty$
and $(\zeta:\psi)=0$ if and only if $\zeta=0$.

Fix once and for all a function $\chi\in\p$ such that $\chi=1$ on $X_0$,
and define a map $L_\psi:\nn\to{\mathbb R}$ by
$$
L_\psi(\zeta)=(\zeta:\psi)/(\chi:\psi).
$$
 
It is routine to show the following properties of $L_\psi$.

\begin{eqnarray}
 L_\psi(c\zeta)=cL_\psi(\zeta),\ \forall c\geq0,\\
 L_\psi(\zeta_1+\zeta_2)\leq L_\psi(\zeta_1)+L_\psi(\zeta_2),\\
  \zeta_1\leq\zeta_2\Rightarrow L_\psi(\zeta_1)\leq L_\psi(\zeta_2),\\
 \s\zeta\subset\D(\gamma)\Rightarrow L_\psi(\zeta\circ\gamma^{-1})
=L_\psi(\zeta),\\
L_\psi(\zeta)\geq 1/(\chi:\zeta).
\end{eqnarray}

\begin{lemma} \label{l1}
If $\eta>0$ and $\xi, \xi'\in\nn$ satisfies $\xi+\xi'=\chi$, then
there is $\delta>0$ such that if $\psi\in\p$, ${\rm diam}(\s\psi)<\delta$
 and
$\zeta\in\nnn$ we have
$$
L_\psi(\xi\zeta)+L_\psi(\xi'\zeta)\leq (1+2\eta)L_\psi(\zeta).
$$
\end{lemma}

\noindent
{\bf Proof.} Given $\eta>0$, there is $\epsilon>0$ such that if
$x, x'\in X_0$ and $d(x,x')\leq \epsilon$, then
$\abs{\xi(x)-\xi(x')}\leq \eta$.
Also this implies $\abs{\xi'(x)-\xi'(x')}\leq \eta$.
 Choose $\delta=\delta(\epsilon)$.  Let
 $\psi$ be as in the lemma and assume
\begin{equation} \label{e1}
\zeta\leq \sum_{i}c_i\psi\circ\gamma_i^{-1}.
\end{equation}

Notice that if we restrict $\gamma_i$ in (\ref{e1}) to 
$\D(\gamma_i)
\cap\s\psi\cap\gamma_i^{-1}(\s\zeta)$, still the inequality
(\ref{e1})
holds. Hence if we choose $x_i$ from $\R(\gamma_i)\subset
\s(\zeta)\subset X_0$,
then for any $x\in\R(\gamma_i)$, we have
$$
\abs{\xi(x)-\xi(x_i)}\leq\eta\ \ {\rm and} \ \ \abs{\xi'(x)-\xi'(x')}\leq \eta.
$$
Moreover the following inequality
$$
\xi(x)\psi\circ\gamma_i^{-1}(x)\leq(\xi(x_i)+\eta)\psi\circ\gamma_i^{-1}(x)
$$
holds for any $x\in X$, since if $x\not\in\R(\gamma_i)$ the both
hand
sides are 0.
Then we have
\begin{eqnarray*}
\zeta(x)\xi(x)\leq\sum_ic_i\xi(x)\psi\circ\gamma_i^{-1}(x)\\
\leq\sum_ic_i(\xi(x_i)+\eta)\psi\circ\gamma_i^{-1}(x).
\end{eqnarray*}
This shows
$$
(\zeta\xi:\psi)\leq\sum_ic_i(\xi(x_i)+\eta).
$$
We have a similar inequality for $\xi'$. Since
$x_i\in X_0$  and thus $\xi(x_i)+\xi'(x_i)=1$, we have
$$
(\zeta\xi:\psi)+(\zeta\xi':\psi)\leq(2\eta+1)\sum_ic_i.
$$
The lemma follows from this.
\hfill
q.\ e.\ d.

\bigskip

Continuing the proof of Theorem \ref{t1}, let us extend the operator
$L_\psi:\nnn\to\mathbb R$ to $C_c(X_0)$ by just putting
$$
L_\psi(\zeta)=L_\psi(\zeta_+)-L_\psi(\zeta_-),
$$
where $\zeta_+$ (resp.\ $\zeta_-$) is the positive (resp.\ negative) part
of $\zeta$.

Then we have:
\begin{equation}
\abs{L_\psi(\zeta)}\leq \Vert \zeta\Vert_\infty, \ \ \forall \zeta\in\nnn.
\end{equation}
In fact if $\zeta\geq 0$, then $\zeta\leq\Vert\zeta\Vert_\infty\chi$,
and thus $L_\psi(\zeta)\leq\Vert\zeta\Vert_\infty$, the general case
following easily from this.

Let us identify $L_\psi$ with the following point of a compact
Hausdorff space:
$$
L_\psi=\{L_\psi(\zeta)\}_\zeta\in\prod_{\zeta\in C_c(X_0)}
[-\Vert\zeta\Vert_\infty,\Vert\zeta\Vert_\infty].
$$

Let $\{\psi_n\}$ be a sequence in $\p$
 such that ${\rm diam}(\s\psi_n)\to0$.
Choose an operator 
$L\in\bigcap_{m}{\rm Cl}\{L_{\psi_n}\mid n\geq m\}$.
This means that for any finite number of elements $\zeta_\nu\in C_c(X_0)$
 and any $\epsilon>0$, there is
a sequence $n_i\to\infty$ such that $\abs{L(\zeta_\nu)-L_{\psi_{n_i}}
(\zeta_\nu)}<\epsilon$. Now we have the following properties of
the map $L:C_c(X_0)\to\mathbb R$.
\begin{eqnarray}
 L(c\zeta)=cL(\zeta),\ \ \forall c\in{\mathbb R},\\
 L(\zeta_1+\zeta_2)\leq L_(\zeta_1)+L(\zeta_2),
\ \ \forall \zeta_1,\ \zeta_2\geq 0, \label{e0}\\
  \zeta_1\leq\zeta_2\Rightarrow L(\zeta_1)\leq L(\zeta_2),\\
 \s\zeta\subset\D(\gamma),\ \gamma\in\Gamma_0\Rightarrow L(\zeta\circ\gamma^{-1})
=L(\zeta),\label{e2}\\
\zeta\in \pp \Rightarrow L(\zeta)\geq 1/(\chi:\zeta), \label{e3}\\
 \abs{L_\psi(\zeta)}\leq \Vert\zeta\Vert_\infty \label{e5}.
\end{eqnarray}

 Moreover by
Lemma \ref{l1} and (\ref{e0}), we have

\begin{lemma} \label{l2}
 If $\zeta\in \nnn$ and $\xi,\ \xi'\in \nn$ satisfy $\xi+\xi'=\chi$,
then
$$
L(\xi\zeta)+L(\xi'\zeta)=L(\zeta).
$$
\end{lemma}

From this one can derive the linearity of $L$. First of all notice that

\begin{equation} \label{e4}
 \zeta,\zeta'\in\nnn\Rightarrow \abs{L(\zeta)-L(\zeta')}\leq\Vert\zeta
-\zeta'\Vert_\infty.
\end{equation}
In fact we have
\begin{eqnarray*}
L(\zeta')=L(\zeta+\zeta'-\zeta)\leq L(\zeta+(\zeta'-\zeta)_+)\leq
L(\zeta)+L((\zeta'-\zeta)_+)\\
\leq L(\zeta)+\Vert (\zeta'-\zeta)_+\Vert_\infty\leq L(\zeta)+\Vert\zeta'-\zeta
\Vert_\infty.
\end{eqnarray*}

Continuing the proof of the linearity, notice that it suffices to show
it
only for those functions $\zeta_1, \zeta_2\in \nnn$.
Choose $\epsilon>0$ small and let
$$\xi_j=(\zeta_j+\epsilon\chi)/(\zeta_1+\zeta_2+2\epsilon)
$$ for $j=1,2$.
Then we have $\xi_1+\xi_2=\chi$. Now
$$
\xi_1(\zeta_1+\zeta_2)-\zeta_1=\epsilon(\zeta_2-\zeta_1)/(\zeta_1
+\zeta_2+\epsilon).
$$
Therefore by (\ref{e4}), we have
$$
\abs{L(\xi_1(\zeta_1+\zeta_2))-L(\zeta_1)}\leq\epsilon.
$$
On the other hand by Lemma \ref{l2}, 
$$
L(\xi_1(\zeta_1+\zeta_2))+L(\xi_2(\zeta_1+\zeta_2)))
=L(\zeta_1+\zeta_2).
$$
Since $\epsilon$  is arbitrary, we have obtained
$$
L(\zeta_1)+L(\zeta_2)=L(\zeta_1)+L(\zeta_2),$$
as is requied.

\bigskip

Now $L$, being a positive operator, corresponds to a Radon measure
$\mu$.
By (\ref{e3}), the measure $\mu$ is nontrivial, and since
(\ref{e5})
implies
$$
\inf\{L(\zeta)\mid \zeta\in \nnn,\ \Vert\zeta\Vert_\infty\leq1\}\leq
1,$$ the measure $\mu$ satisfies $\mu(X_0)\leq 1$.
Finally (\ref{e2}) means the $\Gamma_0$-invariance of $\mu$.

\section{The uniqueness}

In this section $\Gamma$ is again an equicontinuous
and minimal pseudo*group of local homeomorphisms
of a locally compact metric space $X$. The modulus of equicontinuity
is also denoted by $\epsilon\to\delta(\epsilon)$. Denote by
$B_r(x)$ the open $r$-ball in $X$ centered at $x\in X$.

We make the following additional assumption on the pseudo*group $\Gamma$.

\begin{assumption} \label{a1}
 There is a relatively compact open subset $X_0$ of $X$ and $a>0$
such that if $\gamma\in\Gamma$, $x\in X_0$, $x\in\D(\gamma)\subset B_a(x)$ 
and $\gamma x\in X_0$, then there is $\hat\gamma\in\Gamma$ such that
$\D(\hat\gamma)=B_a(x)$ and $\hat\gamma\vert_{\D(\gamma)}=\gamma$.
\end{assumption}

The purpose of this section is to show the following theorem.

\begin{theorem} \label{t3}
Let $\Gamma$ be an equicontinuous and minimal pseudo*group on $X$
satisfying
Assumption \ref{a1}. Then the $\Gamma$-invariant Radon measure on
$X$ is unique up to a scaling.
\end{theorem}

First of all let us show that the holonomy pseudo*group $\Gamma$ on a cross
section $X$ of a minimal
lamination on a compact space $Z$, equicontinuous w.\ r.\ t.\ $X$ 
satisfies Assumption \ref{a1}.
Choose any relatively compact open subset $X_0$ of $X$. 

On one hand by compactness of $Z$ there is $L>0$ such that the germ of any element
of the restriction $\Gamma_0$ to $X_0$ 
is a finite composite of the holonomy maps along leaf curves
of length $\leq$ $L$ that join two points in $X_0$.
 On the other hand there is $a'>0$ such that each leaf
curve of length $\leq L$
starting at $x\in X_0$ and ending at a point in $ X_0$ 
admits a holonomy map defined on the ball $B_{a'}(x)$.
An easy induction shows that Assumption \ref{a1} is satisfied
for $a=\delta(a')$.

\bigskip

Let us embark upon the proof of Theorem \ref{t3}. 
Choosing $a$ even smaller, one may assume that there is
a nonempty open subset $X_1$ of $X_0$ such that the $a$-neighbourhood
$B_a(x)$ of any point $x$ of $X_1$ is contained in $X_0$ and that if
$\gamma\in\Gamma$ and $x'\in X_0$ satisfies  
$\D(\gamma)=B_a(x')$ and $\gamma x'\in X_1$,
then the image $\R(\gamma)=\gamma(B_a(x'))$ is contained in $X_0$.
Choose $b>0$ so that $b\leq \delta(a/3)$, and assume 
there is $x_0\in X_1$ such that $C={\rm Cl}(B)\subset X_1$,
where $B=B_b(x_0)$.

Let $M$ be the space of continuous maps from $C$ to $X_0$,
with the supremum distance $d_\infty$. Define
$$
\Gamma_C=\{\gamma\vert_C \mid \gamma\in \Gamma, \  C\subset\D(\gamma),\
\gamma C\subset X_0\}$$
and let $G$ be the closure of $\Gamma_C$ in $M$.

\begin{lemma} \label{l3}
(1) $G$ is a locally compact metric space.
\\
(2) Any $g\in G$ is a homeomorphism onto a compact subset $gC$ in $X_0$ and
$g$, as well as the inverse map $g^{-1}$, is $\delta(\epsilon)$-continuous.
\end{lemma}

\noindent
{\bf Proof.} All that needs proof is the $\delta(\epsilon)$-continuity 
of $g^{-1}$.
Assume  $\gamma_n\in\Gamma_C$ converge to $g\in G$ in the $d_\infty$-distance.
If $x,x'\in C$ satisfy $d(x,x')>\epsilon$, then $d(\gamma_nx,
\gamma_nx')\geq\delta(\epsilon)$ by the equicontinuity of the
inverse map $\gamma_n^{-1}$. Thus $d(gx,gx')\geq \delta(\epsilon)$,
as is required.
\hfill q.\ e.\ d.

\bigskip

Recall the notations $B=B_b(x_0)$ and $C={\rm Cl}(B)$.

\newcommand{\Cl}{{\rm Cl}}

\begin{lemma} \label{l5}
If $g_n\to g$ in $G$, and $y\in gB$, then for any large $n$
we have $y\in g_nB$ and $g_n^{-1}y\to g^{-1}y$.
\end{lemma}

\noindent
{\bf Proof.}
Choose an arbitrary point $x\in B$ and $\epsilon>0$ such that
${\rm Cl}(B_\epsilon(x))\subset B$.
First let us show that for any $\gamma\in \Gamma_C$,

\begin{equation} \label{e11}
B_{\delta(\epsilon)}(\gamma x)\subset \gamma 
{\rm Cl}(B_\epsilon(x)).
\end{equation}

In fact, by the choice of the number $b$, we have
 $\gamma(B)\subset {\rm Cl}(B_{a/3}(\gamma x_0))$.
That is, $\gamma(B)\subset B_a(\gamma x)$, and thus 
$(\gamma\vert_B)^{-1}$ admits an  extension $\widehat{\gamma^{-1}}\in\Gamma$
defined on $B_a(\gamma x)$. Choose an arbitrary point 
$y\in B_{\delta(\epsilon)}(\gamma x)$. Then by the $\delta(\epsilon)$-continuity
of $\widehat{\gamma ^{-1}}$, the point $x'=\widehat{\gamma^{-1}}y$
lies in ${\rm Cl}(B_{\epsilon}(x))\subset B$. On the other hand 
$x'=\gamma^{-1}\gamma x'=\widehat{\gamma^{-1}}\gamma x'$. Since 
$\widehat{\gamma^{-1}}$
is injective, we have $y=\gamma x'$. This finishes the proof of
 (\ref{e11}).

\bigskip

Next let us show that for any $g\in G$, we have

\begin{equation} \label{e12}
B_{\delta(\epsilon)/2}(gx)\subset g{\rm Cl}(B_\epsilon(x)).
\end{equation}
Again assume  $\gamma_n\in\Gamma_C$ converge to $g\in G$.
Since $\gamma_nx\to gx$, we have for any large $n$ that
$B_{\delta(\epsilon)/2}(gx)\subset B_{\delta(\epsilon)}(\gamma_nx)$.
Thus if $y\in B_{\delta(\epsilon)/2}(gx)$, then by (\ref{e11})
$y=\gamma_nx_n$ for some 
$x_n\in {\rm Cl}(B_\epsilon(x))$. Passing to
a subsequence, assume that $x_n\to x'\in {\rm Cl}(B_\epsilon(x))$.
Now in the following inequality
$$
d(gx',y)=d(gx',\gamma_nx_n)\leq d(gx',\gamma_nx')
+d(\gamma_nx',\gamma_nx_n),
$$
both terms of the RHS can be arbitrarily small if $n$ is sufficiently
large.
That is, $y=gx'$, showing (\ref{e12}).

\bigskip

To finish the proof of the lemma, assume $g_n\to g \in G$ and $y\in gB$.
By (\ref{e12}), for any sufficiently small $\epsilon>0$
we have 
$B_{\delta(\epsilon)/2}(g_ng^{-1}y)\subset g_n{\rm
Cl}(B_\epsilon(g^{-1}y))$.
Since $g_ng^{-1}y\to y$, 
 we have $y\in g_n{\rm Cl}(B_\epsilon(g^{-1}y))$ for any large $n$
and therefore $g_n^{-1}y\in {\rm Cl}(B_\epsilon(g^{-1}y))$.
Since $\epsilon$ is arbitrarily small, this shows the lemma. 
\hfill q.\ e.\ d.

\bigskip

Let $\Gamma_0$ be the restriction of the psudogroup $\Gamma$ to $X_0$.
We shall construct
a pseudo*group $\Gamma_\sharp$ of local homeomorphisms of $G$.
For any $\gamma\in\Gamma_0$, define
\begin{eqnarray*}
& \D(\gamma_\sharp)=\{ g\in G\vert gC \subset \D(\gamma)\}, \\
& \R(\gamma_\sharp)=\{ g\in G\vert gC \subset \R(\gamma)\}, \\
&  \gamma_\sharp g=\gamma\circ g,\ \  \forall g\in \D(\gamma_\sharp.)
\end{eqnarray*}

It may happen that for some $\gamma\in\Gamma_0$,
$\D(\gamma)=\R(\gamma)=\emptyset$. In that case $\gamma_\sharp$
is not defined.

\begin{lemma} \label{l6}
The subsets $\D(\gamma_\sharp)$ and $\R(\gamma_\sharp)$ are open in
$G$, and $\gamma_\sharp$ is $\delta(\epsilon)$-continuous w.\ r.\ 
t.\ the metric $d_\infty$.
\end{lemma}

{\bf Proof.} The easy proof is omitted.
\hfill q.\ e.\ d.

\bigskip

Denote by $\Gamma_\sharp$ the pseudo*group consisting 
of all the
elements $\gamma_\sharp$ for $\gamma\in\Gamma_0$ and their
restrictions to open subsets of the domains.

\begin{lemma} \label{l7}
The action of $\Gamma_\sharp$ on $G$ is minimal.
\end{lemma}

{\bf Proof.}
First let us show that for $\gamma_1,\gamma_2\in\Gamma_C$, there is
$\gamma_\sharp\in\Gamma_\sharp$ such that $\gamma_1\in\D(\gamma_\sharp)$
and that $\gamma_\sharp(\gamma_1)=\gamma_2$.
Since $\gamma_1C\subset B_a(\gamma_1x_0)$, there is an element 
$\gamma'\in\Gamma$ defined on $B_a(\gamma_1x_0)$ which extends
$\gamma_2\circ\gamma_1^{-1}$. Let $\gamma\in\Gamma_0$ be the restriction
of $\gamma'$ to $\Gamma_0$, i.\ e.\ the restriction such taht 
$\D(\gamma)=B_a(\gamma_1x_0)\cap
X_0\cap \gamma'^{-1}X_0$. Clearly $\gamma_1C$ is contained in
$\D(\gamma)$, showing the claim.

Thus we have shown that $\Gamma_\sharp$-orbit of $id_C$ is nothing but
$\Gamma_C$ and hence dense in $G$. To finish the proof, we shall show
that for any $g\in G$, the $\Gamma_\sharp$-orbit of $g$ visits
an arbitrarily small neighbourhood of any element $\gamma_2\in\Gamma_C$.
Let $\epsilon$ be any small number such that the
$2\epsilon$-neighbourhood
of $\gamma_2C$ is contained in $X_0$. Take $\gamma_1\in\Gamma_C$ such
that
$d_\infty(g,\gamma_1)<\delta(\epsilon)$. Choosing $\epsilon$ and hence
$\delta(\epsilon)$ even smaller, one may very well assume that $gC$
is contained in $B_a(\gamma_1x_0)$.
Then the element
$\gamma\in\Gamma_0$ constructed above (for $\gamma_1$ and $\gamma_2$)
contains $gC$ in its domain, i.\ e.\ $g$ is contained in
$\D(\gamma_\sharp)$, and furthermore $d_\infty(\gamma_\sharp g,
\gamma_2)<\epsilon$.
\hfill q.\ e.\ d.

\bigskip

Now by Lemmata \ref{l3}, \ref{l6} and \ref{l7}, one can apply
Theorem \ref{t2} to $(\Gamma_\sharp, G)$ to find a
$\Gamma_\sharp$-invariant
Radon measure $m$ on $G$. 
(This is the point where the concept of pseudo*group is useful.
Notice that even if $\Gamma_\sharp$ is equicontinuous,
it does not necessarily imply that
the pseudogroup generated by $\Gamma_\sharp$ is equicontinuous.)
One can assume $m$ is a probability measure
since $G$ is in fact a precompact open subset of a bigger space.
Now let $\mu$ and $\mu'$ be distinct $\Gamma_0$-invariant
probability measures on $X_0$. Then their restrictions to $B$ are
also distinct, by the minimality of the $\Gamma_0$-action.
That is, there is a function $\zeta\in C_c(B)_{}$ such that 
$\mu(\zeta)\neq\mu'(\zeta)$. One may assume further that
$\zeta$ is nonnegative valued.

\begin{lemma} \label{l8}
For any $g\in G$, we have
$$
\int_{X_0}\zeta (g^{-1}x)\mu(dx)=\int_{X_0}\zeta(x)\mu(dx).
$$
\end{lemma}

{\bf Proof.} For $g\in\Gamma_C$, this is just the $\Gamma_0$-invariance
of $\mu$. For general $g$, assume $\gamma_n\to g$ for 
$\gamma_n\in \Gamma_C$. Then by Lemma \ref{l5}, if $x\in gB$, then
$x\in\gamma_n B$ for any large $n$ and $\gamma_n^{-1}x\to g^{-1}x$.
If $x\not\in gB$, then since $\gamma_n\s(\zeta)\to g\s(\zeta)$
in the Hausdorff distance, $\zeta(\gamma_n^{-1}x)=0$ for any large $n$,
as well as $\zeta(g^{-1}x)$. In any case for any $x\in X_0$, we have
$\zeta(\gamma_n^{-1}x)\to\zeta(g^{-1}x)$. The lemma follows from
the dominated convergence theorem.
\hfill q.\ e.\ d.

\bigskip

Now recall the space $X_1$. It is an open subset of $X_0$ which
contains $C$ such that
the $a$-neighbourhood
$B_a(x)$ of any point $x$ of $X_1$ is contained in $X_0$ 
and that if
$\gamma\in\Gamma$ and $x'\in X_0$ satisfies  
$\D(\gamma)=B_a(x')$ and $\gamma x'\in X_1$,
then the image $\R(\gamma)=\gamma(B_a(x'))$ is contained in $X_0$.

\begin{lemma} \label{l9}
The function 
$$
Z(x)=\int_G\zeta(g^{-1}x)m(dg)
$$
is constant on $X_1$.
\end{lemma}

{\bf Proof.} Define a function $\zeta_x:G\to {\mathbb R}$ by
$\zeta_x(g)=\zeta(g^{-1}x)$. Lemma \ref{l5} and an additional argument
as above shows that
$\zeta_x$ is a continuous function.
 
Choose $x,x'\in X_1$ on the same $\Gamma$-orbit.
By the assumption of $X_1$, there is $\gamma\in\Gamma_0$
such that $\gamma x=x'$ and $\D(\gamma)=B_a(x)\subset X_0$ and 
$\R(\gamma)\subset X_0$. Then
we have
$$
\{g\in G\mid\zeta_x(g)>0\}\subset \D(\gamma_\sharp).
$$
In fact if $\zeta_x(g)=\zeta(g^{-1}x)>0$, then $x\in gB$.
On the other hand, diam$(gB)\leq 2a/3$, and thus 
$gC\subset B_a(x)=\D(\gamma)$, i.\ e.\ $g\in\D(\gamma_\sharp)$.

By the $\Gamma_\sharp$-invariance of the measure $m$, we have
\begin{eqnarray*}
Z(x)=\int_G\zeta_x(g)m(dg)=\int_G\zeta_x(\gamma_\sharp^{-1}(g))m(dg)
=\int_G\zeta_x(\gamma^{-1}\circ g)m(dg)\\
=\int_G\zeta(g^{-1}\gamma x)m(dg)
=\int_g\zeta_{\gamma x}(g)m(dg)=Z(x')
\end{eqnarray*}

That is, the function $Z$ is constant along a $\Gamma$-orbit in
$X_1$. On the other hand it is continuous, since $\zeta\circ g^{-1}$
has the same modulus of continuity. Now the minimality of
$\Gamma_0$-action on $X_1$ shows the lemma.
\hfill q.\ e.\ d.

\bigskip

\begin{lemma} \label{l10}
The function $Z$ is constant on $X_0$.
\end{lemma}

{\bf Proof.} It suffices to show that for
 any $x'\in X_0$ and $x\in X_1$ on the same $\Gamma_0$-orbit, we have
$Z(x)=Z(x')$.  By the assumption of $X_1$, there exists
an element $\gamma\in\Gamma_0$ such that $\gamma x'= x$ and 
$\D(\gamma)=B_a(x')\cap X_0$. Then just as before, one can show
$$
\{g\in G\mid\zeta_{x'}(g)>0\}\subset \D(\gamma_\sharp).
$$
Again by the $\Gamma_\sharp$-invariance of $\mu$, we have
$Z(x)=Z(x')$.
\hfill q.\ e.\ d. 

\bigskip

Now let us finish the proof of Theorem \ref{t3}. By Lemma \ref{l10},
the function $Z$ is  constant on $X_0$, depending only on $\zeta$
and $m$.
We have on one hand
$$
\int_{X_0}\int_G\zeta(g^{-1}x)m(dg)\mu(dx)=\int_{X_0}Z\mu(dx)=Z.
$$
On the other hand by Fubini and by Lemma \ref{l8}
$$
Z=\int_G\int_{X_0}\zeta(g^{-1}x)\mu(dx)m(dg)=\int_G\mu(\zeta)m(dg)=\mu(\zeta).
$$

Since $Z$ does not depend on the choice of $\mu$, we have
$\mu(\zeta)=\mu'(\zeta)$, contrary to the assumption.

\section{The uniqueness of harmonic measures for group actions}

Here the notations of the previous sections are all abandoned.
Let $\alpha:\Gamma\times X\to X$ be an effective (i.\ e.\ faithful) action of 
a countable group $\Gamma$ on a
compact metric space $X$, and let $p$ be a probability measure
on $\Gamma$, i.\ e.\ a function $p:\Gamma\to[0,1]$
such that $\sum_{\gamma\in\Gamma}p(\gamma)=1$. We assume
that $\s(p)=\{\gamma\in\Gamma\mid p(\gamma)>0\}$ 
generates $\Gamma$ as a semigroup.
A probability measure $\mu$ on $X$ is called {\em $p$-harmonic}
if $\mu=\alpha_*(p\times\mu)$, that is, for any continuous function
$f$ on $X$, we have
$$
\int_Xf(x)\mu(dx)=\int_X\sum_{\gamma\in\Gamma}p(\gamma)f(\gamma x)\mu(dx).
$$

This section is devoted to the proof of the following theorem.

\begin{theorem} \label{t4}
If the action $\alpha$ is equicontinuous and minimal, then
the $p$-harmonic probability measure $\mu$ on $X$ is unique.
\end{theorem}

{\bf Proof.}
Let $M$ be the space of continuous maps from $X$ to $X$,
endowed with the supremum metric $d_\infty$, and let $G$ be the closure
of $\Gamma$ in $M$.
Then as in section 3, Lemmata \ref{l3} and \ref{l5}, 
we can show that $G$ is a compact metrizable
topological group, (with the topology induced from the metric
$d_\infty$).

Let $f$ be an arbitrary continuous function on $X$.
Let $m$ be a Haar probability measure on $G$. Define a function 
$f_m:X\to {\mathbb R}$
by
$$
f_m(x)=\int_Gf(gx)m(dg).
$$
The function $f_m$ is on one hand continuous since the functions $f\circ g$
have the same modulus of continuity, and on the other hand constant
on $\Gamma$-orbits by the right invariance of $m$. Hence by the
minimality
of the action, $f_m$ is
a constant, which we denote by $c(f, m)$.

Let $\mu$ be a $p$-harmonic probability measure on $X$, and define a function
$f_\mu:G\to {\mathbb R}$ by
$$
f_\mu(g)=\int_Xf(gx)\mu(dx).
$$
Then $f_\mu$ is a continuous function w.\ r.\ t.\ $d_\infty$,
and by the $p$-harmonicity of $\mu$, it satisfies
$$
f_\mu(g)=\sum_{\gamma\in\Gamma}p(\gamma)f_\mu(g\gamma ).
$$
If $f_\mu$ takes the maximal value at $g\in G$, then for
any $\gamma\in\s(p)$, the value of $f_\mu$ at $g\gamma $ is also
the maximal. Repeating this arguments one can show that
$f_\mu$ takes the maximal value on the coset $g\Gamma $, since
$\s(p)$ generates $\Gamma$ as a semigroup. 
Because $\Gamma$ is dense in $G$, the function $f_\mu$ must be
a constant, equal to $f_\mu(e)=\mu(f)$.

Now we have
$$
\int_X\int_Gf(gx)m(dg)\mu(dx)=\int_X c(f,m) \mu(dx)=c(f,m).$$
On the other hand by Fubini,
$$
c(f,m)=\int_G\int_Xf(gx)\mu(dx)m(dg)=\int_G\mu(f)m(dg)=\mu(f).$$

Then the value $\mu(f)$, being equal to $c(f,m)$, does not
depened on $\mu$, showing the uniqueness of $\mu$.
\hfill q.\ e.\ d.


\bigskip

\begin{thebibliography}{99}

\bibitem[C]{C} A. Candel, {\em The harmonic measures of Lucy Garnett,}
Adv.\ Math.\ {\bf 176}(2003) 187-247.

\bibitem[DKN]{DKN} B. Deroin, V. Kleptsyn and A. Navas {\em
Sur la dynamique unidimensionelle en r\'egularit\'e interm\'ediaire,}
Acta Math.\ {\bf 199}(2007) 199-262.

\bibitem[G]{G} L. Garnett, {\em Foliations, the ergodic theorem and
Brownian motion,} J. Funct.\ Anal.\ {\bf 51}(1983) 285-311.

\bibitem[S]{S} R. Sacksteder, {\em Foliations and pseudogroups,}
Amer.\ J. Math.\ {\bf 87}(1965) 79-102.

\end{thebibliography}
\end{document}